\documentclass[10pt]{amsart}

\newtheorem{thm}{Theorem}[section]

\newtheorem{obs}{Remark}[section]
\newtheorem{defin}{Definition}[section]

\usepackage{amssymb}
\usepackage{amsmath}
\usepackage{amsfonts}
\usepackage{graphicx}

\numberwithin{equation}{section}

\usepackage[T1]{fontenc}
\usepackage{palatino}

\begin{document}
\title[vortex stretching and anisotropic diffusion]
{Vortex stretching and anisotropic diffusion\\ 
in the 3D Navier-Stokes equations}
\author{Z. Gruji\'c}
\address{Department of Mathematics\\
University of Virginia\\ Charlottesville, VA 22904}
\dedicatory{Dedicated to Professor Hugo Beir\~ao da Veiga on the occasion
     of his $70$th birthday, with admiration.}
\date{\today}

\begin{abstract}

The goal of this article is to present -- in a cohesive, and somewhat self-contained 
fashion -- several
recent results revealing an experimentally, numerically, and mathematical
analysis-supported \emph{geometric scenario} manifesting \emph{large data} logarithmic 
\emph{sub-criticality} of the 3D Navier-Stokes regularity problem.
Shortly -- in this scenario --
the \emph{transversal small scales} produced by the mechanism of vortex stretching 
(coupled with the decay of the volume of the regions of intense vorticity)
reach the threshold sufficient for the \emph{locally anisotropic diffusion} to engage and 
control the sup-norm of the vorticity, preventing the (possible) formation of 
finite time singularities.

\end{abstract}

\maketitle

\vspace{.1in}

\section{Prologue}

\noindent Vortex stretching has been viewed as the principal \emph{physical mechanism}
responsible for the vigorous creation of \emph{small scales} in turbulent fluid flows. This
goes back at least to G. I. Taylor's fundamental paper ``Production and dissipation of vorticity in
a turbulent fluid'' from 1937 \cite{Tay37}. 

\medskip

While the \emph{production} part has been relatively well-understood (the amplification of
the vorticity via the process of vortex stretching follows essentially from the conservation of
the angular momentum in the incompressible fluid), the precise physics/mathematics behind
the vortex stretching-induced \emph{dissipation} is less transparent.
For his part, Taylor inferred the thoughts on the \emph{anisotropic dissipation} chiefly
from the wind tunnel \emph{measurements} of turbulent flow past a uniform grid, concluding
the paper with the following sentence.

\medskip

``It seems that the stretching of vortex filaments must be regarded as the principal
mechanical cause of the high rate of dissipation which is associated with
turbulent motion.''

\medskip

Since then, it has been a grand challenge in the mathematical fluid mechanics community to
try to explain/quantify the process of anisotropic dissipation in turbulent flows directly from the 
mathematical model -- the 3D Navier-Stokes equations (NSE).

\medskip

Numerical simulations (cf. \cite{AKKG87, JWSR93, SJO91, VM94}) reveal that the regions of 
intense vorticity 
are dominated by \emph{coherent vortex structures} and in particular, \emph{vortex
filaments}. There are two imminent morphological signatures of this geometry. One
is  \emph{local coherence of the vorticity direction}, and the other one is 
\emph{local existence of sparse/thin direction(s)}. 

\medskip

\noindent \texttt{Local coherence. Geometric depletion of the nonlinearity.}  

\medskip

\noindent The pioneering 
work in this direction was presented by 
Constantin in \cite{Co94} where he
derived a singular integral representation of the stretching factor in the evolution of the
vorticity magnitude featuring a geometric kernel that is depleted by \emph{local coherence}
of the \emph{vorticity direction}, a purely geometric condition. This has been referred to as
\emph{geometric depletion of the nonlinearity}, and has led to the first rigorous confirmation
of the phenomenon of anisotropic dissipation utilizing the 3D NSE, a theorem (\cite{CoFe93}) stating
that as long as the vorticity direction is Lipschitz-coherent (in the regions of high
vorticity), the $L^2$-norm of the vorticity is controlled, and no finite time blow-up can
occur. The Lipschitz-coherence condition was later scaled down to $\frac{1}{2}$-H\"older
in \cite{daVeigaBe02}, and a full \emph{spatiotemporal localization} of the $\frac{1}{2}$-H\"older condition
was performed in \cite{Gr09} (a different approach to localization was previously introduced in
\cite{ChKaLe07}).
A family of local, hybrid geometric-analytic regularity criteria including a
\emph{scaling invariant} improvement of the $\frac{1}{2}$-H\"older condition was presented
in \cite{GrGu10-1}. 
The study of the coherence of the vorticity direction up to the boundary-regularity
criteria in the case of the no-stress boundary conditions was presented in
\cite{daVeigaBe09}, and in the case of the no-slip boundary conditions in \cite{daVeiga07}.

\medskip

Essentially, an unhappy event preventing the manifestation of the anisotropic
dissipation in this setting is the one of `crossing of the vortex lines', i.e., of the vorticity direction 
forming even a simple spatial discontinuity 
-- two different limit points --
at a (possible) singular time. 

\medskip
 
\noindent \texttt{Local anisotropic sparseness. Vortex stretching-anisotropic diffusion.} 

\medskip

\noindent An alternative mathematical description of
the anisotropic dissipation in the 3D incompressible
viscous flows was recently exposed in \cite{Gr13}, and is based on the concept of \emph{anisotropic
diffusion}. Taken at face value, the 3D NSE diffusion -- generated by the Laplacian --
is isotropic. 
The (isotropic) diffusion is then utilized via sharp, local-in-time \emph{spatial analyticity} properties of
solutions in $L^\infty$ which provide an ambient amenable to the application of the
\emph{harmonic measure majorization principle}. More precisely,
as long as the region of intense vorticity,
defined to be the region in which the vorticity magnitude exceeds a fraction of the $L^\infty$
norm, exhibits the property of \emph{local existence of a sparse/thin direction} at a scale comparable to
the radius of spatial analyticity (essentially, $\frac{1}{C} \|\omega(t)\|_\infty^{-\frac{1}{2}}$
where
$\omega$ denotes the vorticity of the fluid),
an argument relying on the translational and rotational invariance
of the equations and certain geometric properties of the harmonic measure
(this is what 
introduces \emph{anisotropy}),
and the harmonic measure maximum principle, shows that the $L^\infty$ 
norm of the vorticity is controlled, and no finite-time blow up can occur. It is worth mentioning
that it suffices to assume the aforementioned sparseness property 
\emph{intermittently in time}. 

\medskip

Of course, a key question is whether there is any evidence, either numerical, or
mathematical, that the \emph{scale of local linear sparseness/thinness} needed for 
triggering the mechanism of anisotropic diffusion is in fact achieved
in a turbulent flow. Thinking in terms of vortex filaments, the scale we are interested in is
essentially the \emph{length scale of the diameters of the cross-sections}, i.e.,
the \emph{transversal scale} of the filament.
It appears easier -- both numerically and analytically -- to estimate the \emph{axial length scale 
of the filaments} instead. This, coupled with a suitable estimate on the volume of the region
of the intense vorticity, provides an (indirect) estimate on the desired length scale.

\medskip

Direct numerical simulations suggest that -- intermittently in time/in the time average --
the axial lengths of the filaments are essentially comparable to the \emph{macro scale}
(e.g., the side length $L$ in the case of the $L$-periodic boundary conditions). 
On the other hand, the \emph{a priori} $L^1$-estimate on the vorticity \cite{Co90} implies that the
volume of the region of intense vorticity is bounded by 
 $C \|\omega(t)\|_\infty^{-1}$. Hence (intermittently in time),
 the transversal \emph{micro scale} of the filament is bounded by
 $C \|\omega(t)\|_\infty^{-\frac{1}{2}}$; in other words, the  NSE regularity problem
 in this scenario becomes \emph{critical}.
 
 \medskip
 
 In addition to the numerical evidence, a very recent work \cite{DaGr12-3} presented a \emph{mathematical
 evidence} of creation and persistence (in the time average) of the macro scale-long 
 vortex filaments. More precisely, exploiting a \emph{dynamic, spatial multi-scale ensemble 
 averaging
 process} designed to detect \emph{sign-fluctuations} of an \emph{a priori}
 sign varying physical quantity across scales, it was shown that there exists a range of scales -- extending
 from a suitable micro scale to the macro scale -- at which the vortex stretching
 term is essentially positive. (The averaging procedure utilized had been previously developed in a
 recent series of papers \cite{DaGr11-1, DaGr11-2, DaGr12-1, DaGr12-2} as a mathematical framework for the study
 of turbulent cascades in \emph{physical scales} of 3D incompressible fluid flows.)

 \medskip
 
 The aforementioned ruminations offer a physically, numerically, and 
 mathematical analysis-supported \emph{large data criticality scenario} for the
 3D NSE. The NSE themselves are (still) \emph{super-critical}; regardless of the
 functional setup, there has been a `scaling gap' between a regularity criterion
 in view and the corresponding \emph{a priori} bound. An instructive example is
 given by the regularity condition $u \in L^\infty_t L^3_x$ obtained by Escauriaza, Seregin
 and Sverak in \cite{ESS03}, to be contrasted to
 Leray's \emph{a priori} bound $u \in L^\infty_t L^2_x$ (\cite{Le34}).
  
 \medskip
 
 A natural question to ask is whether it is possible to \emph{break the criticality}
in this setting; i.e., whether the intricate interplay between the vortex stretching and the anisotropic
diffusion results in preventing the formation of singularities, rather than in a critical blow-up
 scenario. In a very recent article \cite{BrGr13-2}, it was shown
 that a \emph{very mild, purely geometric} assumption yields a
uniform-in-time $L \log L$ bound on the vorticity; this in turn implies an extra-log decay
of the vorticity distribution function, i.e., of the volume of the region of intense
vorticity, breaking the scaling, and transforming the aforementioned criticality scenario
into an \emph{anisotropic diffusion-win scenario} (no singularities).
More precisely, the assumption is a
uniform-in-time boundedness of the localized \emph{vorticity direction} in a suitable, logarithmically
weighted, local space of \emph{bounded mean oscillations} ($BMO$).
An interesting feature of this space 
is that it allows for
discontinuous functions exhibiting singularities of, e.g.,  $\sin \log |\log ( \, \mbox{something algebraic} \, )|$-type.
Hence, the vorticity direction can form a singularity in a geometrically spectacular fashion
-- every point on the unit sphere being a limit point -- and the $L \log L$ bound will still hold
(in particular, a simple `crossing of the vortex lines' is not an obstruction).

\medskip

The proof is based
on an adaptation of the method utilized in \cite{Co90}, the novel components being exploiting
\emph{analytic cancelations}
in the vortex-stretching term via a version of the Div-Curl Lemma (in the sense of
Coifman, Lions, Meyer and Semmes theory of compensated compactness in Hardy spaces), a local
version of the $\mathcal{H}^1-BMO$ duality, a sharp pointwise multiplier theorem in local $BMO$ spaces, and 
Coifman-Rochberg's  $BMO$-estimate on the logarithm of the \emph{maximal function} of a locally 
integrable function (the estimate is \emph{independent of the function} and depending only 
on the dimension of the space).

\medskip

This result (\cite{BrGr13-2}) is -- in a way -- complementary to the results obtained in \cite{BrGr13-1}. 
The class of conditions leading to an $L \log L$-bound presented in \cite{BrGr13-1} 
consists of suitable \emph{blow-up} rates that
can be characterized
as `wild in time' with a uniform spatial (e.g., algebraic) structure, while the condition presented
in \cite{BrGr13-2} can be characterized as `wild in space' and uniform in time.

\medskip

In summary, the papers \cite{Gr13, DaGr12-3, BrGr13-2} can be viewed as providing a rigorous mathematical
framework (directly from the 3D NSE) for justification of Taylor's view on vortex stretching 
as the principal
mechanical cause for the high rate of dissipation in turbulent flows. 
Incidentally, they also point to a possible new
direction in the study of the 3D NSE regularity problem.

\section{Anisotropic diffusion}

\noindent 3D Navier-Stokes equations (NSE) -- describing a flow of 3D
incompressible viscous fluid -- read

\[
 u_t+(u\cdot \nabla)u=-\nabla p + \nu \triangle u,
\]
supplemented with the incompressibility condition $ \, \mbox{div} \,
u = 0$, where $u$ is the velocity of the fluid, $p$ is the
pressure, and $\nu$ is the viscosity. Taking the curl yields the vorticity formulation,

\[
 \omega_t+(u\cdot \nabla)\omega= (\omega \cdot \nabla)u + \nu \triangle
 \omega,
\]
where $\omega = \, \mbox{curl} \, u$ is the vorticity.
Utilization of the identity $\triangle u = - \, \mbox{curl} \, \omega$ leads
to the Biot-Savart law,
\[
 u(x) = c \int \nabla\frac{1}{|x-y|} \times \omega(y) \, dy,
\]
closing the system for the vorticity field.

\medskip

Computational simulations of \emph{3D homogeneous turbulence}
reveal that the regions of intense
vorticity organize in \emph{coherent vortex structures}, and in particular,
in elongated vortex tubes/filaments, cf.
\cite{S81, AKKG87, SJO91, JWSR93, VM94}.
An in-depth analysis of creation and dynamics of vortex tubes in 3D
turbulent flows was presented in \cite{CPS95}; 
in particular, a suitably defined \emph{dynamical scale of coherence}
of the vorticity direction field was estimated. The current body of work containing
analytical, as well as analytical and numerical results on the dynamics
of coherent vortex structures includes \cite{GGH97, GFD99, Oh09, Hou09}.

\medskip

In what follows, we will focus on  \emph{sparseness}.

\medskip

\begin{defin}\emph{
Let $x_0$ be a point in $\mathbb{R}^3$, $r>0$, $S$ an open subset of
$\mathbb{R}^3$ and $\delta$ in $(0,1)$.
The set $S$ is \emph{linearly $\delta$-sparse around $x_0$ at scale
$r$ in weak sense} if there exists a unit vector $d$ in $S^2$ such
that
\[
 \frac{|S \cap (x_0-rd, x_0+rd)|}{2r} \le \delta.
\]
}
\end{defin}

Denote by $\Omega_t(M)$ the vorticity super-level set at time $t$;
more precisely,
\[
 \Omega_t(M) = \{x \in \mathbb{R}^3: |\omega(x,t)| > M\}.
\]
Then the following manifestation of anisotropic diffusion
holds (\cite{Gr13}).

\begin{thm}\label{sparse_omega}
Suppose that a solution $u$ is regular on an interval $(0,T^*)$.

\medskip

Assume that either

\medskip

(i) \ there exists $t$ in $(0,T^*)$ such that
$\displaystyle{t+\frac{1}{d_0^2 \|\omega(t)\|_\infty} \ge T^*}$, 
or

\medskip

(ii) \ $\displaystyle{t+\frac{1}{d_0^2 \|\omega(t)\|_\infty} < T^*}$
for all $t$ in $(0,T^*)$,
and there exists $\epsilon$ in $(0,T^*)$ such that for any
$t$ in $(T^*-\epsilon, T^*)$, there
exists $s=s(t)$ in $\Bigl[t+\frac{1}{4d_0^2 \|\omega(t)\|_\infty},
t+\frac{1}{d_0^2 \|\omega(t)\|_\infty}\Bigr]$
with the property that for any spatial point
$x_0$, there exists a scale $r=r(x_0)$, $0<r\le \frac{1}{2d_0^2
\|\omega(t)\|_\infty^\frac{1}{2}}$, such that the super-level set
$\Omega_s(M)$ is linearly $\delta$-sparse around $x_0$ at scale $r$
in weak sense; here, $\delta=\delta(x_0)$ is an arbitrary value
in $(0,1)$,
$h=h(\delta)=\frac{2}{\pi}\arcsin\frac{1-\delta^2}{1+\delta^2}$,
$\alpha=\alpha(\delta)\ge\frac{1-h}{h}$, and
$M=M(\delta)=\frac{1}{d_0^\alpha}
\|\omega(t)\|_\infty$.

\medskip

Then, there exists $\gamma >0$ such that $\omega$ is in
$L^\infty\Bigl((T^*-\epsilon, T^*+\gamma); L^\infty\Bigr)$, i.e.,
$T^*$ is not a singular time. \ \ ($d_0$ is a suitable absolute constant.)
\end{thm}

The quantity $(d_0^2 \|\omega(t)\|_\infty)^{-1}$ is the time step in 
the local-in-time well-posedness scheme in $L^\infty$ initiated
at $t$. The scheme can be complexified \cite{GrKu98, Gu10}. For any
$s \in \bigl(t, t+(d_0^2 \|\omega(t)\|_\infty)^{-1}\bigr)$, $\omega(s)$
is a restriction of the function holomorphic in the region
$\{x+iy \in \mathbb{C}^3: \, |y| < 1/c_1 \sqrt{s}\}$; moreover, the
sup-norm of the complexified solution is controlled by $c_2 \|\omega(t)\|_\infty$.

\medskip

The idea of the proof is as follows. Let $x_0$ be a spatial point, and $d=d(x_0)$
a sparse direction within the region of intense vorticity, at the scale comparable to 
the uniform lower bound on the
radius of spatial analyticity. By the translational invariance of the equations,
we can send $x_0$ to the origin, and by the rotational invariance, we can
align $d$ with one of the coordinate directions. The (real)
coordinate is then embedded in the complex plane, and the 
\emph{harmonic measure maximum principle} applied with respect to the disk centered at the
origin -- with the radius comparable to the analyticity radius -- is utilized to
exploit the sparseness condition resulting in a ``self-improving'' bound on the sup-norm
of the complexified vorticity, preventing the finite time blow-up.

\medskip

The main engine behind the argument is local-in-time spatially analytic smoothing in $L^\infty$,
a strong manifestation of the (isotropic) diffusion generated by $\partial_t - \triangle$; a \emph{locally
anisotropic diffusion effect} is a consequence of the 
translational and rotational invariance of the equations, and
geometric properties of the harmonic measure.

\begin{obs}\emph{
It suffices to assume the sparseness condition at (suitably chosen)
\emph{finitely many} times/intermittently in time.}
\end{obs}

\section{A possible road to criticality.}

\noindent  Adopting the notation introduced in the preceding section, define \emph{the region of intense vorticity} to
be the set
$\displaystyle{
\Omega_{s(t)}\Bigl(\frac{1}{c_1} \|\omega(t)\|_\infty\Bigr)}$
for an appropriate $c_1>1$. Let $R_0$ be a suitable \emph{macro scale} associated with the flow.
Computational simulations indicate that (intermittently-in-time) dominant geometry in
the region of intense vorticity is the one of $R_0$-long vortex filaments; in
order to estimate the transversal micro-scale of the filaments, it suffices to
have a good estimate on the rate of the decrease of the volume of the vorticity
super-level sets.

\medskip

Let $(0,T)$ be an interval of interest.
In \cite{Co90},
provided the initial vorticity is a bounded measure (and the initial velocity
is of finite energy),
Constantin showed that a corresponding weak solution satisfies
$\displaystyle{\sup_{t \in (0,T)} \|\omega(t)\|_{L^1} \le c_{0,T} = c(u_0, \omega_0, T)}$.
Chebyshev's inequality then implies
\[
 \, \mbox{Vol} \, \biggl( \Omega_{s(t)} \Bigl(\frac{1}{c_1}
\|\omega(t)\|_\infty \Bigr) \biggr) \le \frac
{c_{0,T}'}{\|\omega(t)\|_\infty} \ \ (c_{0,T}'>1),
\]
which -- in turn -- yields the decrease of the transversal micro-scale
of the filaments of at least
at least $\displaystyle{\frac
{c_{0,T}''}{\|\omega(t)\|^\frac{1}{2}_\infty}}$ $(c_{0,T}''>1)$.
This is precisely the scale of local, linear sparseness needed to trigger
the mechanism of anisotropic diffusion exposed in the previous 
section, i.e., we arrive at \emph{criticality}.

\medskip

It is instructive to check the scaling in the \emph{geometrically worst case
scenario} -- no sparseness -- the super level set being clumped in a ball.
In this case, the criticality would require

\[
 \lambda_{\omega(t)}(\beta) = O \biggl(\frac{1}{\beta^{3/2}}\biggr)
\]
uniformly in $(T^*-\epsilon, T^*)$ ($\lambda$ denotes the distribution function);
this is a scaling-invariant condition -- back to super-criticality of the problem,

\[
 O \biggl(\frac{1}{\beta^{3/2}}\biggr) \ \ {vs.} \ \  O \biggl(\frac{1}{\beta^1}\biggr).
\]
(In fact, this is precisely the vorticity analogue of 
the velocity scaling gap -- $L^\infty_t L^3_x$ \emph{vs.} 
$L^\infty_t L^2_x$.)

\medskip

Summarizing -- in this scenario -- the vortex stretching acts as the mechanism
bridging (literally) the scaling gap in the regularity problem.

\section{Mathematical evidence of criticality.}

\noindent In this section, we identify the
range of scales of \emph{positivity of the vortex-stretching term}
$S \omega \cdot \omega$ ($S$ denotes the symmetric part of
the gradient of $u$); this corresponds to the range of scales of
creation and persistence of vortex filaments.

\medskip

To this end, we exploit a spatial \emph{multi-scale averaging method} designed to
detect \emph{sign fluctuations} of a quantity of interest across 
\emph{physical scales} recently introduced in the study of turbulent transport rates 
in 3D incompressible fluid flows \cite{DaGr11-1, DaGr11-2, DaGr12-1, DaGr12-2}.

\medskip

Let $B(0,R_0)$ be a macro-scale domain.
A \emph{physical scale} $R$, $0 < R \le R_0$, is realized via suitable 
ensemble averaging of the
localized quantities with respect to `$(K_1,K_2)$-covers at scale
$R$'.

\begin{defin}\emph{
Let $K_1$ and $K_2$ be two positive integers, and $0 < R \le R_0$; a
cover $\{B(x_i,R)\}_{i=1}^n$ of $B(0,R_0)$ is a
\emph{$(K_1,K_2)$-cover at scale $R$} if
\[
 \biggl(\frac{R_0}{R}\biggr)^3 \le n \le K_1
 \biggr(\frac{R_0}{R}\biggr)^3,
\]
and any point $x$ in $B(0,R_0)$ is covered by at most $K_2$ balls
$B(x_i,2R)$.
}
\end{defin}

The parameters $K_1$ and $K_2$ are the maximal \emph{global}
and \emph{local multiplicities}, respectively.

\medskip

For a physical density of interest $f$, consider -- 
localized to the cover elements $B(x_i,
R)$ (per unit mass) -- local quantities $\hat{f}_{x_i,R}$,
\[
\hat{f}_{x_i,R} =  \frac{1}{R^3}
\int_{B(x_i,2R)} f(x) \psi^\delta_{x_i,R} (x) \, dx,
\]
for some $0 < \delta \le 1$. The smooth cut-off
functions $\psi_i=\psi_{x_i,R}$ are equal to 1 on $B(x_i,R)$, vanish outside
of $B(x_i, 2R)$, and satisfy
\[
 |\nabla \psi_i| \le c_\rho \frac{1}{R} \psi_i^\rho, \ \ \ 
 |\triangle \psi_i| \le c_\rho \frac{1}{R^2} \psi_i^{2\rho-1},
\]
for a suitably chosen $\rho$,  $\frac{1}{2} < \rho < 1$. Denote by $\psi_0$
the cut-off corresponding to the macro-scale domain $B(0,R_0)$.
The cut-offs associated with `boundary elements', i.e., the cover
elements $B(x_i,R)$ intersecting the boundary of the macro-scale 
domain, are modified to satisfy certain compatibility relations with
the global cut-off $\psi_0$; for technical details see, e.g., \cite{DaGr11-1}.

\medskip

Denote by $\langle F\rangle_R$
the \emph{ensemble average} given by
\[
 \langle F\rangle_R = \frac{1}{n} \sum_{i=1}^n
 \hat{f}_{x_i,R}.
\]
The key feature of  $\{\langle
F\rangle_R\}_{0<R\le R_0}$ is that $\langle F\rangle_R$ being
\emph{stable} -- i.e., nearly-independent on a particular choice of
the cover (with the fixed local multiplicity $K_2$) -- indicates
there are \emph{no significant sign fluctuations} at scales comparable or
greater than $R$. On the other hand, if $f$ does exhibit significant sign fluctuations
at scales comparable or greater than $R$, suitable
\emph{rearrangements} of
the cover elements up to the maximal multiplicity -- emphasizing
first the positive and then the negative parts of $f$ --
will result in $\langle F\rangle_R$ experiencing a wide range of
values, from positive through zero to negative, respectively
(the larger $K_2$, the finer detection).

\medskip

For a non-negative density $f$, the ensemble
averages are all
comparable to each other throughout the full range of scales, $0
< R \le R_0$; in particular, they are all comparable to the simple
average over the macro-scale domain,

\begin{equation}
 \frac{1}{K_1} F_0 \le \langle F \rangle_R \le K_2 F_0,
\end{equation}
for all $0 < R \le R_0$, where
\[
 F_0= \frac{1}{R_0^3} \int  f(x)
 \psi_0^\delta (x) \, dx.
\]

\medskip

Back to vortex stretching.
Denote the time-averaged localized vortex-stretching terms per unit
mass associated to the cover element $B(x_i,R)$ by $VST_{x_i,R,t}$,

\begin{equation}
VST_{x_i,R,t} = \frac{1}{t} \int_0^t \frac{1}{R^3} \int (\omega
\cdot \nabla)u \cdot  \omega \ \phi_i \, dx \, ds;
\end{equation}
here, $\phi_i=\phi_i(x,s)=\eta(s) \, \psi_{x_i,R}(x)$, where $\eta$ is a smooth function
on the time-interval of interest $(0,T)$ satisfying 
\[
 \eta=0 \ \mbox{on} \ (0, 1/3 \, T), \ \ \eta=1 \ \mbox{on} \ (2/3 \, T, T), \ \
 |\eta'| \le c_\kappa \frac{1}{T} \eta^\kappa,
\]
for a suitable $\kappa$, $0 < \kappa < 1$.
The quantity of interest is the ensemble average of
$\{VST_{x_i,R,t}\}_{i=1}^n$,
\begin{equation}\label{PhiR}
 \langle VST\rangle_{R,t} = \frac{1}{n}\sum_{i=1}^n VST_{x_i,R,t}.
\end{equation}

\medskip

$B(x_i,R)$-localized \emph{enstrophy level dynamics} is then as follows,

\begin{align}\label{loc}
 \int_0^t \int (\omega \cdot \nabla)u \cdot \phi_i \, \omega \; dx
 \; ds
 &=
 \int \frac{1}{2}|\omega(x,t)|^2\psi_i(x) \; dx + \int_0^t \int
 |\nabla\omega|^2\phi_i
 \;  dx \; ds\notag\\
 &- \int_0^t \int \frac{1}{2}|\omega|^2 \bigl((\phi_i)_s+\triangle\phi_i\bigr) \; dx
 \; ds\notag\\
 &- \int_0^t \int \frac{1}{2}|\omega|^2 (u \cdot \nabla\phi_i) \; dx
 \; ds,
\end{align}
for any $t$ in $(2/3 \, T, T)$, and $1 \le i \le n$. The hope is that the
ensemble-averaging of the right-hand side 
detects a \emph{dynamic range of scales} of 
positivity of the vortex stretching term $S \omega \cdot \omega$.

\medskip

Before stating the result, several macro-scale quantities need to
be introduced. Denote by $E_{0,t}$ time-averaged enstrophy per unit 
mass associated
with the macro scale domain $B(0,2R_0) \times (0,t)$,
\[
 E_{0,t}=\frac{1}{t}\int_0^t \frac{1}{R_0^3} \int \frac{1}{2}|\omega|^2
 \phi_0^{1/2} \, dx \, ds,
\]
by $P_{0,t}$ a modified time-averaged palinstrophy per unit mass,
\[
 P_{0,t}= \frac{1}{t}\int_0^t \frac{1}{R_0^3} \int |\nabla\omega|^2
 \phi_0 \, dx \, ds
 + \frac{1}{t}\frac{1}{R_0^3} \int \frac{1}{2}|\omega(x,t)|^2
 \psi_0(x) \, dx
\]
(the modification is due to the shape of the temporal cut-off
$\eta$), and by $\sigma_{0,t}$ a corresponding Kraichnan-type scale,
\[
 \sigma_{0,t}=\biggl(\frac{E_{0,t}}{P_{0,t}}\biggr)^\frac{1}{2}.
\]
Then the following holds \cite{DaGr12-3}.

\begin{thm}
Let $u$ be a global-in-time local Leray solution on $\mathbb{R}^3
\times (0,\infty)$, regular on $(0,T)$. Suppose that, for some $t
\in (2/3 \, T, T)$,
\begin{equation}\label{cond}
C \max\{M_0^\frac{1}{2}, R_0^\frac{1}{2}\} \,
\sigma_{0,t}^\frac{1}{2} < R_0
\end{equation}
where $\displaystyle{M_0=\sup_t \int_{B(0,2R_0)} |u|^2 < \infty}$,
and $C > 1$ a suitable constant depending only on the cover
parameters.

Then,
\begin{equation}\label{estimate}
 \frac{1}{C} \, P_{0,t} \le \langle VST \rangle_{R,t} \le
 C \, P_{0,t}
\end{equation}
for all $R$ satisfying
\begin{equation}\label{range}
C \max\{M_0^\frac{1}{2}, R_0^\frac{1}{2}\} \,
\sigma_{0,t}^\frac{1}{2} \le R \le R_0.
\end{equation}
\end{thm}

\medskip

\noindent \texttt{A couple of remarks.}

\medskip

\noindent (i) Suppose that $T$ is the first (possible) singular time, and that the
macro-scale domain contains some of the spatial singularities (at
time $T$). This, paired with the assumption that $u$ is a
global-in-time local Leray solution implies
\[
 \sigma_{0,t} \to 0, \  t \to T^- ;
\]
hence, the condition (\ref{cond}) in the theorem is \emph{automatically
satisfied} for any $t$ sufficiently close to the singular time $T$.

\medskip

\noindent  (ii) $P_{0,t} \to \infty, \  t \to T^-$, i.e., the vortex stretching intensifies 
as we approach the singularity.

\medskip

\noindent (iii) The power of $\frac{1}{2}$ on $\sigma_{0,t}$ is a correction originating in 
the need for a suitable control of the localized 
transport term.

\section{Logarithmic sub-criticality}

\noindent The purpose of this section is to show how a \emph{very mild, purely geometric},
condition transforms the criticality scenario exposed in the previous sections
into a log sub-critical scenario, preventing the possible formation of singularities.

\medskip

The idea is to get a uniform-in-time $L \log L$-bound on $w$,
where $w=\sqrt{1+|\omega|^2}$; this would impose an extra decay on the distribution 
function of the vorticity, breaking the criticality.

\medskip

Suppose that the solution in view is smooth on $(0,T)$, and focus on some
macro-scale spatial domain, e.g., $B(0,R_0)$. 
The evolution of $w$ 
satisfies the following partial differential inequality (\cite{Co90}),

\begin{equation}\label{w}
\partial_t w - \triangle w + (u \cdot \nabla) w \le \omega \cdot \nabla u \cdot 
\frac{\omega}{w}.
\end{equation}
Since our goal is to control the evolution of $w \log w$ over $B(0,R_0)$, it is convenient to
multiply (\ref{w}) by $\psi \, (1+\log w)$ where $\psi=\psi_0$ (a smooth cut-off associated
with the macro-scale domain as in the previous section).
After a fair amount of calculation, the
following bound transpires,
\[
\begin{aligned}
I(\tau) \equiv \int \psi(x) \, w(x,\tau) \log w(x,\tau) \, dx &\le 
 I(0) + c \int_0^\tau \int_x \omega \cdot \nabla u \cdot \psi \, \xi  \, \log w \, dx \, dt\\
 &+ \ \mbox{\emph{a priori} \, bounded},
\end{aligned}
\]
for any $\tau$ in $[0,T); \, \xi$ denotes the vorticity direction.

\medskip

Note that with respect to the scaling of \emph{a priori} bounded quantities, the integral on 
the right-hand side is log super-critical.
The strategy to overcome this is as follows. Exploit the
\emph{analytic cancellations} in the vortex-stretching term $\omega \cdot \nabla u$
utilizing a version of the Div-Curl Lemma in the Hardy space $\mathcal{H}^1$, and
then transform the gain via $\mathcal{H}^1-BMO$ duality into some sort of an
oscillation condition on $\xi$. To do this in an efficient manner, 
we will need a sharp pointwise multiplier theorem in a version of local $BMO$, 
and a result quantifying the intimate relationship between $\log$ and $BMO$.

\medskip

In order to keep the exposition self-contained, a number of relevant definitions 
and results from harmonic analysis are listed below; for more details, as well as
the references, see \cite{BrGr13-2}.

\medskip

\noindent Let $f$ be a distribution. The maximal function of $f$ is defined as
\[
M_h f(x)=\sup_{t>0}|f*h_t(x)|,  \ x \in \mathbb{R}^n,
\] 
where $h$ is a fixed, normalized test function supported in the unit ball, and 
$h_t$ denotes $t^{-n}h(\cdot/t)$. A distribution $f$ is in the Hardy space 
$\mathcal {H}^1$ if $\|f\|_{\mathcal {H}^1}=\|M_h f\|_1<\infty$.
The local maximal function is defined as,
\[
m_h f(x)=\sup_{0<t<1}|f*h_t(x)|, \ x \in \mathbb{R}^n;
\]
a distribution $f$ is in the local Hardy space $\frak{h}^1$ if $\|f\|_{\frak{h}^1}=\|m_h f\|_1
<\infty$.

\medskip

\noindent Div-Curl Lemma \ Suppose that
$E$ and $B$ are $L^2$-vector fields satisfying $\mbox{div} \,  E = \, \mbox{curl} \,  B
= 0$  (in the sense of distributions).  Then,
\[
\|E\cdot B\|_{\mathcal H^1} \leq c(n) \, \|E\|_{L^2} \|B\|_{L^2}.
\]

\medskip

\noindent The space of bounded mean oscillations, $BMO$, is defined as follows
\[
 BMO = \biggl\{ f \in L^1_{loc} : \, \sup_{x \in \mathbb{R}^n, r>0}
 \Omega \bigl(f, I(x,r)\bigr) < \infty \biggr\}
\]
where $\displaystyle{\Omega \bigl(f, I(x,r)\bigr)=\frac{1}{|I(x,r)|}\int_{I(x,r)} |f(x)-f_I | \, dx}$
is the mean oscillation of the function $f$ with respect to its mean
$f_I = \frac{1}{|I(x,r)|}\int_{I(x,r)} f(x) \, dx$, over the cube $I(x,r)$ centered at $x$ with
the side-length $r$. A local version of $BMO$, usually denoted by $bmo$, is defined by 
finiteness of the following
expression,
\[
 \|f\|_{bmo} = \sup_{x \in \mathbb{R}^n, 0 < r < \delta} \Omega \bigl(f, I(x,r)\bigr)
 + \sup_{x \in \mathbb{R}^n, r \ge \delta} \frac{1}{|I(x,r)|} \int_{I(x,r)} |f(y)| \, dy,
\]
for some positive $\delta$.

\medskip

\noindent $\bigl(\mathcal{H}^1\bigr)^* = BMO$ and $\bigl(\frak{h}^1\bigr)^* = bmo$; 
the duality is realized via integration.

\medskip

\noindent When $f \in L^1$, we can focus on small scales, e.g., $0 < r < \frac{1}{2}$.
Let $\phi$ be a positive, non-decreasing function on $(0, \frac{1}{2})$, 
and consider the following version of local
weighted spaces of bounded mean oscillations,
\[
 \|f\|_{\widetilde{bmo}_\phi} = \|f\|_{L^1} +  \sup_{x \in \mathbb{R}^n, 0<r<\frac{1}{2}}
 \frac{\Omega \bigl(f, I(x,r)\bigr)}{\phi(r)}.
\]
Of special interest will be the spaces $\widetilde{bmo} = \widetilde{bmo}_1$, and 
$\widetilde{bmo}_{\frac{1}{|\log r|}}$

\medskip

\noindent Let $h$ be in $\widetilde{bmo}$, and $g$ in 
$L^\infty \cap \widetilde{bmo}_\frac{1}{|\log r|}$. Then,
\[
 \|g \, h\|_{\widetilde{bmo}} \le c(n) \, \Bigl( \|g\|_\infty + \|g\|_{\widetilde{bmo}_\frac{1}{|\log r|}} \Bigr) 
 \ \|h\|_{\widetilde{bmo}}.
\]

\medskip

\noindent Let $M$ denote
the Hardy-Littlewood maximal operator, and $f$ be a locally integrable function. Then,
\[
 \| \log M f \|_{BMO} \le c(n)
\]
(the bound is completely
\emph{independent of} $f$). The estimate remains valid if we replace $M f$ with 
$\mathcal{M} f = \bigl( M \sqrt{|f|}\bigr)^2$;
the advantage of working with $\mathcal{M}$ is that the $L^2$-maximal theorem implies
the following estimate
\[
 \|\mathcal{M} f\|_1 \le c(n) \|f\|_1.
\]

\medskip

Back to $J \equiv \displaystyle{\int_0^\tau \int_x \omega \cdot 
\nabla u \cdot \psi \, \xi  \, \log w \, dx \, dt}$. Decomposing $\log w$ as
\[
\log w = \log \frac{w}{\mathcal{M}w} 
+ \log \mathcal{M}w,
\]
induces a decomposition of the integral; denote this decomposition by 
$J=J_1+J_2$.

\medskip

The first integral can be bounded without paying much attention to oscillations, 
utilizing a couple of elementary inequalities and the $L^2$-maximal 
theorem.

\medskip

For the second integral, we have the following string of inequalities,

\[
\begin{aligned}
J_2
 & = \int_0^\tau \int_x \omega \cdot \nabla u \cdot \psi \, \xi \,
  \log \mathcal{M}w \, dx \, dt\\
 & \le c \int_0^\tau \|\omega \cdot \nabla u\|_{\frak{h}^1} 
 \|\psi \, \xi \log \mathcal{M}w\|_{bmo} \, dt\\
 & \le c \int_0^\tau \|\omega \cdot \nabla u\|_{\mathcal{H}^1} 
 \|\psi \, \xi \log \mathcal{M}w\|_{\widetilde{bmo}} \, dt\\
  & \le c \int_0^\tau \|\omega\|_2 \|\nabla u\|_2 \Bigl( \|\psi \, \xi\|_\infty + \|\psi \, \xi\|_{\widetilde{bmo}_{\frac{1}{|\log r|}}}
  \Bigr)
 \Bigl(\|\log \mathcal{M}w\|_{BMO} + \|\log \mathcal{M}w\|_1\Bigr) \, dt\\
 & \le c \sup_{t \in (0,T)} \ \biggl\{ \Bigl(1 + \|\psi \, \xi\|_{\widetilde{bmo}_{\frac{1}{|\log r|}}}\Bigr)  
  \Bigl(\|\log \mathcal{M}w\|_{BMO} + \|\log \mathcal{M}w\|_1\Bigr)
     \biggr\} \ \ \int_t \int_x |\nabla u|^2 \\ 
 &\le c \ \Bigl(1 + \sup_{t \in (0,T)} \|\psi \, \xi\|_{\widetilde{bmo}_{\frac{1}{|\log r|}}}\Bigr) \                               
 \Bigl(1+ \sup_{t \in (0,T)} \|\omega\|_1\Bigr) \  
 \ \int_t \int_x |\nabla u|^2.
\end{aligned}
\]

\medskip

Since the $L^1$-norm of the vorticity is \emph{a priori} bounded, we arrive at the
following theorem \cite{BrGr13-2}.

\begin{thm}
Let $u$ be a Leray solution to the 3D NSE. Assume that the initial vorticity $\omega_0$ 
is in $L^1 \cap L^2$, and that $T>0$ is the first (possible) blow-up time. Suppose that
\[
 \sup_{t \in (0,T)} \| (\psi \xi) (\cdot, t)  \|_{\widetilde{bmo}_\frac{1}{|\log r|}} < \infty.
\]
Then,
\[
 \sup_{t \in (0,T)} \int \psi(x) \, w(x,t) \log w(x,t) \, dx < \infty.
\]
\end{thm}

\medskip

\noindent \texttt{Good news.} $\widetilde{bmo}_\phi$ contains discontinuous
functions if and only if
$\displaystyle{\int_0^\frac{1}{2} \frac{\phi(r)}{r} \, dr = \infty}$. 
More specifically, 
$\widetilde{bmo}_{\frac{1}{|\log r|}}$ contains bounded functions
with the discontinuities of, say, 
$\displaystyle{\sin \log |\log (\, \mbox{something algebraic} \, )|}$-type, i.e., 
$\xi$ can (as it approaches $T$, and the spatial singularity at $T$)
oscillate among \emph{infinitely many limit points} on the unit sphere,
and still yield extra-log decay of the distribution function of $\omega$
breaking the criticality.

\section{Epilogue}

\noindent From the fluid mechanics perspective, the results reviewed provide
a framework for \emph{rigorous identification} of the interplay between vortex
stretching and anisotropic diffusion as a principal mechanism behind
the phenomenon of \emph{turbulent dissipation}.

\medskip

From the PDE perspective, they identify a \emph{large data geometric sub-criticality
scenario} in the 3D NS regularity problem. This is achieved in two steps. First, a 
\emph{dynamic criticality scenario} is 
revealed -- thinking in terms
of vortex filaments -- in which the transversal scale of the filaments matches the scale
of local, linear (anisotropic) sparseness of the region of intense vorticity needed to trigger 
the anisotropic diffusion \cite{Gr13, DaGr12-3} . Then,
a very mild geometric condition -- boundedness of the vorticity direction in space
$\displaystyle{\widetilde{bmo}_{\frac{1}{|\log r|}}}$ -- \emph{breaking the criticality} is identified. 
In particular, the vorticity direction is allowed to \emph{develop spatial discontinuities}
at the possible singular time $T$ \cite{BrGr13-2}.
It is instructive to briefly compare this to several (relatively) recent results from the literature in 
which a form of
\emph{criticality} is \emph{assumed}, and then 
an \emph{anisotropic condition} implying the regularity is identified.
In \cite{SeSv09} (see also \cite{CSYT08, CSTY09}), it is shown -- under the
type I blow-up assumption -- that the local \emph{axisymmetric solutions} do not form
singularities. The regularity condition here is the one of \emph{global anisotropy}.
In \cite{GiMi11}, the authors showed -- also under the type I blow-up assumption --
that as long as the vorticity direction possesses a uniform
\emph{modulus of continuity}, no finite time blow-up can occur. As in
\cite{BrGr13-2}, the regularity condition here is the one of
\emph{local anisotropy}; however, in contrast to \cite{BrGr13-2}, 
uniform continuity of the vorticity direction is still required.

\medskip

\bigskip

\bigskip

\noindent ACKNOWLEDGMENTS \ The author acknowledges support of the \emph{National Science
Foundation} via the grant DMS-1211413 and the \emph{Research Council of Norway} via
the grant F20/213473.

\end{document}